\documentclass[12pt]{amsart}
\usepackage{amsfonts}
\usepackage{amssymb}

\theoremstyle{definition}

\theoremstyle{remark}

\newcommand{\ds}{\displaystyle}

\begin{document}

\title[On Bochner flat almost K\"ahler manifolds ]
{On Bochner flat almost K\"ahler manifolds}%

\thanks
{This research is supported by the National Science Fund 
of the Bulgarian Ministry of Education and Science 
under Grant~DFNI-T01/0001, 2012.}

\author{Ognian Kassabov}%
\address{University of Transport, Sofia, Bulgaria}
\email{okassabov@abv.bg}

\begin{abstract}
The main purpose of this article is to prove that there exist no proper $AK_3$-manifold
of dimension $2n$, $n\ge 3$, with vanishing Tricerri-Vanhecke Bochner curvature tensor
and constant scalar curvature. 
\end{abstract}

\maketitle

\section{Introduction}

The Riemannian manifolds, that are conformal to a flat Riemannian manifold, form
a very important class. As it is well known, a characteristic property
of these manifolds (in dimension $\ge 3$) is that they have vanishing Weil 
curvature tensor.

Studying the Betti numbers of K\"ahler manifolds, Bochner defined for them 
a tensor, as an algebraic analogue of 
the Weil tensor. Although we don't yet know the exact geometric meaning of 
the Bochner curvature tensor, it is an object of special interest, because in many 
cases in the K\"ahler geometry (not only in the study of Betti numbers) this tensor 
plays a role similar to that of the Weil tensor for Riemannian manifolds. 
For example a Riemannian manifold is of constant sectional curvature
if and only if it is Einsteinian and has vanishing Weil curvature tensor.
Analogously a K\"ahler manifold is of constant holomorphic sectional curvature if and 
only if it is Einsteinian and has vanishing Bochner curvature tensor.   

In 1981 Tricerri and Vanhecke defined a Bochner-type curvature tensor for an
arbitrary almost Hermitian manifold. In particular, in the K\"ahler case
their Bochner tensor coincides with the classical Bochner one. Since then, many studies 
have been also made about this tensor for different classes of almost Hermitian manifolds.

K\"ahler manifolds with vanishing Bochner curvature tensor and constant scalar
curvature are classified in \cite{M-T}. The same problem for nearly K\"ahler 
manifolds is studied in\cite{OK1}, see also  \cite{EPS}. In the present paper we consider 
the case of almost K\"ahler manifolds satisfying the third curvature condition. Namely, in 
section 4 we prove that such a manifolds of dimension $2n$, $n\ge3$, must be K\"ahler (Theorem 3). 

In section 3 we prove for semi-K\"ahler manifolds a result, similar to the well known theorem
of Tricerri and Vanhecke for the so-called generalized complex space forms \cite{T-V}.




\setcounter{equation}{0}
\section{Preliminaries}

In this section $x,y,z,u,v$ will be arbitrary vectors in a point $p$ of a $2n$-dimensional 
almost Hermitian $(AH)$  manifold $ M$ with metric tensor $g$  and almost complex structure $J$. 
The curvature tensor (of type (1,3) or ((0,4)), the Ricci tensor (of type (1,1) or 
((0,2)) and the scalar curvature are denoted by $R$, $S$ and $\tau$, respectively.

We will use the second Bianchi identity
\begin{equation} \label{eq:2.6}
	(\nabla_xR)(y,z,u,v)+(\nabla_yR)(z,x,u,v)+(\nabla_zR)(x,y,u,v)=0  
\end{equation}
as well as the Ricci  identities
$$
	(\nabla_x(\nabla_yJ))(z)-(\nabla_y(\nabla_xJ))(z)=R(x,y)Jz-JR(x,y)z
$$
$$
	(\nabla_x(\nabla_yP))(z,u)-(\nabla_y(\nabla_xP))(z,u)=-P(R(x,y)z,u)-P(z,R(x,y)u)
$$
where $\nabla$ is the covariant differentiation with respect to the Riemannian connection of $M$
and $P$ is a tensor field of type (0,2). 
 
In the almost Hermitian geometry it is convenient to use the following operators:
$$
	\begin{array}{rl}
    	\varphi(P)(x,y,z,u)&  =g(x,u)P(y,z)-g(x,z)P(y,u)\\
												 & \  +g(y,z)P(x,u)-g(y,u)P(x,z) \ ,
	\end{array}											 
$$
$$
	\begin{array}{rl}
    	\psi(P)(x,y,z,u)&  =g(x,Ju)P(y,Jz)-g(x,Jz)P(y,Ju)-2g(x,Jy)P(z,Ju)\\
												 & \  +g(y,Jz)P(x,Ju)-g(y,Ju)P(x,Jz)-2g(z,Ju)P(x,Jy)
	\end{array}											 
$$
for a tensor field $P$ of type (0,2). Put also $\pi_1=1/2\varphi(g)$, $\pi_2=1/2\psi(g)$. 

Recall  that  the manifold is 
conformal flat if and only if its Weil tensor $C(R)$ vanishes, where
$$
	C(R)=R-\frac{1}{2n-2}\varphi(S)+\frac{\tau}{(2n-1)(2n-2)}\pi_1 \ .
$$

The classes of K\"ahler $(K)$, nearly K\"ahler $(NK)$, almost K\"ahler ($AK$),
quasi K\"ahler $(QK)$, semi-K\"ahler $(SK)$ manifolds are defined respectively 
by $\nabla J=0$, $(\nabla_xJ)x=0$,
\begin{equation} \label{eq:2.7}
	dF=0 \qquad {\rm or} \qquad g((\nabla_xJ)y,z)+g((\nabla_yJ)z,x)+g((\nabla_zJ)x,y)=0 \ ,
\end{equation}
$(\nabla_xJ)y+(\nabla_{Jx}J)Jy=0$, $\delta F=0$, where $F(x,y)=g(x,Jy)$ is the fundamental form
and $\delta $ denotes the coderivative. The following inclusions are strict
$$
	K\subset NK\subset QK\subset SK \hspace{1in}  K\subset AK\subset QK\subset SK \ ,
$$
see e.g. \cite{G-H}. Moreover $K=NK\cap AK$.

For a class $L$ of almost Hermitian manifolds its subclass $L_i$ is defined by  the
$i$-th of the following identities for its curvature tensor
 
1) $R(x,y,z,u)$=$R(x,y,Jz,Ju)$;

2) $R(x,y,z,u)=R(x,y,Jz,Ju)+R(x,Jy,z,Ju)+R(Jx,y,z,Ju)$;

3) $R(x,y,z,u)$=$R(Jx,Jy,Jz,Ju)$.

\noindent 
Then $AH_1\subset AH_2\subset AH_3$, $K=K_1$, $NK=NK_2$.

For almost Hermitian manifolds a second Ricci tensor and a second scalar curvature are introduced.
Namely, the $*$-Ricci tensor $S^*$ and the $*$-scalar curvature $\tau^*$ are given by
$$
	S^*(x,y)=\sum_{i=1}^{2n}\, R(x,e_i,Je_i,Jy) \qquad \tau^*(x,y)=\sum_{i=1}^{2n}\, S^*(e_i,e_i) 
$$ 
where $\{ e_i, i=1,...,2n \}$ is an orthonormal basis of $T_pM$. 
Then for $AH_3$-manifolds the Tricerri-Vanhecke Bochner curvature tensor $B(R)$ of $M$ is defined by
$$
	B(R)=R-(\varphi+\psi)(T)-\varphi(Q)+\mu(\pi_1+\pi_2)  +\nu \pi_1\ ,  
$$
where 
$$
	T=\frac1{8(n+2)}(S+3S^*)-\frac1{8(n-2)}(S-S^*) \qquad\qquad
	Q=\frac1{2(n-2)}(S-S^*)
$$
$$
	\mu=\frac{\tau+3\tau^*}{16(n+1)(n+2)}-\frac{\tau-\tau^*}{16(n-1)(n-2)} \qquad\qquad
	\nu=\frac{\tau-\tau^*}{4(n-1)(n-2)} \ ,
$$
see \cite{T-V}. In particular, if the manifold is K\"ahler, this tensor coincides with 
the classical Bochner curvature tensor. Note that for any $AH_3$-manifold the tensors 
$S$, $S^*$ (and so also $T$ and $Q$) are symmetric and $J$-invariant (hybrid). Note also
that for a manifold $M\in AH_3$ with $B(R)=0$ it follows $M\in AH_2$.

For an $AK_2$-manifold the following identities hold:
\begin{equation} \label{eq:2.2}
	2(\nabla_xQ)(y,z)=Q((\nabla_xJ)y,Jz)+Q(Jy,(\nabla_xJ)z) \ ,
\end{equation}
\begin{equation} \label{eq:2.3}
	R(x,y,z,u)-R(x,y,Jz,Ju)=\frac12g(K(x,y),K(z,u)) \ ,
\end{equation}
where $K(x,y)=(\nabla_xJ)y-(\nabla_yJ)x$, see \cite{Barros}, \cite{Gray}. 
Note that  (\ref{eq:2.2}) and $M\in QK$ imply
\begin{equation} \label{eq:2.4}
	(\nabla_xQ)(y,z)+(\nabla_xQ)(Jy,Jz) =0\ , \qquad (\nabla_xQ)(y,z)+(\nabla_{Jx}Q)(y,Jz) =0\ , 
\end{equation}
\begin{equation} \label{eq:2.5}
	(\nabla_xQ)(y,Jz) =(\nabla_xQ)(Jy,z)=(\nabla_{Jx}Q)(y,z) \ . 
\end{equation}
Moreover, if the scalar curvature $\tau$ of $M$ is a constant, then by (\ref{eq:2.4}) $\tau^*$
is also constant. So $\mu$ and $\nu$ are constants, too.




\setcounter{equation}{0}
\section{A Theorem for semi-K\"ahler manifolds.}

\vspace{0.2cm}
{\bf Theorem 1.} {\it Let $M$ be a $2n$-dimensional semi-K\"ahler manifold, $n>2$,
whose curvature tensor has the form
\begin{equation} \label{eq:4.1}
	R=\varphi(P)+f\pi_1+h\pi_2 \ ,
\end{equation}
where $f$ and $h$ are functions and $P$ is a symmetric tensor field of type (0,2). Then
$h$ is a constant. If $h=0$, then $M$ is conformal flat. 
If $h\ne 0$, then $M$ is a K\"ahler manifold
of constant holomorphic sectional curvature. }

\begin{proof}
First of all we note that (\ref{eq:4.1}) implies (with a contraction) that the tensor field $P$ is symmetric.
Suppose that  $x,y,z$ are  unit vectors in $T_pM$ such that $x,y,z, Jx,Jy,Jz$ are
mutually orthogonal. 

From the second Bianchi identity 
$$
	(\nabla_{Jx}R)(y,Jy,Jz,z)+(\nabla_{y}R)(Jy,x,Jz,z)+(\nabla_{Jy}R)(Jx,y,Jz,z)=0
$$
it follows
\begin{equation} \label{eq:4.2}
	 Jx(h)-h\Big( g((\nabla_yJ)y,x)+g((\nabla_{Jy}J)Jy,x) \Big) =0 \ .
\end{equation}
Since $M$ is semi-K\"ahler, this implies easily  that $h$ is a constant.
If $h=0$ by a standard way we obtain from (\ref{eq:4.1})
$$
	R=\frac1{n-2}\varphi(S)-\frac{\tau}{2(n-1)(2n-2)}\pi_1 \ ,
$$
so $M$ is conformal flat. 

Let $h\ne 0$. Then (\ref{eq:4.2}) becomes 
$$
	 g((\nabla_yJ)y,x)+g((\nabla_{Jy}J)Jy,x) =0 \ .
$$

Now we change in (\ref{eq:2.6}) $(x,y,z,u,v)$ with $(Jx,y,z,Jz,Jy)$ and because of $h=const.\ne 0$ we find
$$
	 g((\nabla_yJ)y,x)+g((\nabla_zJ)z,x) =0 \ .
$$
From the last two equalities  we find $ (\nabla_yJ)y=0$, so
$M$ is a nearly K\"ahler manifold. 

We replace in (\ref{eq:2.6}) $(z,u,v)$ by $(Jz,Jz,z)$. The result is
$$
	(\nabla_xP)(y,z)-(\nabla_yP)(x,z)+
	h\Big( 3g((\nabla_xJ)y,Jz)+3g((\nabla_yJ)Jz,x)+2g((\nabla_{Jz}J)x,y) \Big)=0  \ .
$$
Making a cyclic sum in the above equality and using that $M$ is nearly K\"ahler we derive
$$
	g((\nabla_xJ)y,z)+g((\nabla_yJ)z,x)+g((\nabla_zJ)x,y)=0 
$$ 
for all $x,y,z$  in $T_pM$ such that $x,y,z, Jx,Jy,Jz$ are
mutually orthogonal. Note that  the last equality is  true also 
when $z=y$ or $z=Jy$. Hence it follows that it is  true for arbitrary
$x,y,z$  in $T_pM$, so  $M$ is also almost K\"ahler. Consequently
$M$ is K\"ahlerian, thus proving the assertion. 
\end{proof}

\vspace{0.2cm}
{\bf Remark.} An almost Hermitian manifold is said to be {\it a generalized complex space form} if
its curvature tensor has the form
$$
	R=f\pi_1+h\pi_2 \ ,
$$
where $f$ and $h$ are functions. For such a manifold Tricerri and Vanhecke \cite{T-V} proved
that it is of constant sectional curvature or a K\"ahler manifold of constant holomorphic
sectional curvature.




\setcounter{equation}{0}
\section{The main result.}

In the following five lemmas we assume that $M$ is an $AK_3$-manifold with vanishing 
Tricerri-Vanhecke Bochner curvature tensor and constant scalar curvature.
As noted in section 2 in this case the manifold is $AK_2$ and of constant
$*$-scalar curvature, so $\mu$ and $\nu$ are  also constants.

We begin  by proving that under the above assumptions  
the tensor $T$ has the property of the type (\ref{eq:2.2}).

\vspace{0.2cm}
{\bf Lemma 1.} {\it The tensor $T$ satisfies
\begin{equation} \label{eq:3.5}
	2(\nabla_xT)(y,z)=T((\nabla_xJ)y,Jz)+T(Jy,(\nabla_xJ)z) 
\end{equation}
for all $x,y,z\in T_pM$.}

\begin{proof}
Let $x,y$ be unit vectors in $T_pM$ with $x\perp y,Jy$. 
Putting in the second Bianchi identity (\ref{eq:2.6}) $z=u=Jy$, $v=y$ 
and using (\ref{eq:2.4}) we obtain 
\begin{equation} \label{eq:3.1}
	(\nabla_xT)(y,y) +(\nabla_xT)(Jy,Jy)=(\nabla_yT)(x,y) +(\nabla_{Jy}T)(x,Jy) \ . 
\end{equation}
Analogously if we put in (\ref{eq:2.6}) $z=Jy$, $u=Jx$, $v=x$, we find
$$
	(\nabla_xT)(x,x) +(\nabla_xT)(Jx,Jx)+(\nabla_xT)(y,y) +(\nabla_xT)(Jy,Jy)=4(\nabla_yT)(x,y) +4(\nabla_{Jy}T)(x,Jy) \ . 
$$
From the last two equalities it follows
\begin{equation} \label{eq:3.2}
	(\nabla_xT)(x,x) +(\nabla_xT)(Jx,Jx)=3(\nabla_xT)(y,y) +3(\nabla_{x}T)(Jy,Jy) \ . 
\end{equation}
Let $\{ e_i,Je_i,\, i=1,...n\} $ be an orthonormal basis of $T_pM$ such that $x=e_1$. We
put in (\ref{eq:3.2}) $y=e_i$ and we add for $i=2,...n$. Then using $x(\tau)=x(\tau^*)=0$ we obtain
\begin{equation} \label{eq:3.3}
	(\nabla_xT)(x,x) +(\nabla_xT)(Jx,Jx)=0 \ . 
\end{equation}
Now (\ref{eq:3.2}) becomes
\begin{equation} \label{eq:3.4}
	(\nabla_xT)(y,y) +(\nabla_{x}T)(Jy,Jy) =0  
\end{equation}
for $x\perp y,Jy$.   By (\ref{eq:3.3}) and (\ref{eq:3.4}) we may conclude that (\ref{eq:3.4})
holds for arbitrary vectors $x,y\in T_pM$. Since $T$ is a symmetric tensor, 
this implies
$$
	(\nabla_xT)(y,z) +(\nabla_{x}T)(Jy,Jz) =0  
$$
for all $x,y,z \in T_pM$. Hence using 
$$
	(\nabla_{x}T)(Jy,Jz)=(\nabla_{x}T)(y,z)-T((\nabla_xJ)y,Jz)-T(Jy,(\nabla_xJ)z
$$ 
we obtain the assertion.
\end{proof}

\vspace{0.2cm}
{\bf Remark.} It follows from Lemma 1 that $T$ satisfies also the analogues of (\ref{eq:2.4})
and (\ref{eq:2.5}).

\vspace{0.2cm}
{\bf Lemma 2.} {\it The tensor $T$ satisfies}
$$
	T(R(x,y)z,u)+T(z,R(x,y)u)+T(R(x,y)Jz,Ju)+T(Jz,R(x,y)Ju)
$$
$$
	=\frac12\Big( T((\nabla_xJ)(\nabla_yJ)z,u)+T(z,(\nabla_xJ)(\nabla_yJ)u)
	-T((\nabla_yJ)(\nabla_xJ)z,u)-T(z,(\nabla_yJ)(\nabla_xJ)u) \Big) \ .
$$

\begin{proof}
Using Lemma 1 we calculate
$$
	(\nabla_x(\nabla_yT)(z,u)=\frac12\Big( 
	T((\nabla_x(\nabla_yJ))z,Ju)+T(Jz,(\nabla_x(\nabla_yJ))u)  \Big)
$$
$$
	+\frac14\Big( T((\nabla_xJ)u,(\nabla_yJ)z) +T((\nabla_xJ)z,(\nabla_yJ)u)-T((\nabla_xJ)(\nabla_yJ)z,u)-T(z,(\nabla_xJ)(\nabla_yJ)u)  \Big) \ .
$$
Hence
$$
	\hspace{-4in} (\nabla_x(\nabla_yT)(z,u)-(\nabla_y(\nabla_xT)(z,u)
$$
$$
	=\frac12\Big( 
	T((\nabla_x(\nabla_yJ))z,Ju)+T(Jz,(\nabla_x(\nabla_yJ))u)
	-T((\nabla_y(\nabla_xJ))z,Ju)-T(Jz,(\nabla_y(\nabla_xJ))u)  \Big)
$$
$$
	+\frac14\Big( T((\nabla_yJ)(\nabla_xJ)z,u)+T(z,(\nabla_yJ)(\nabla_xJ)u)-T((\nabla_xJ)(\nabla_yJ)z,u)-T(z,(\nabla_xJ)(\nabla_yJ)u)  \Big) \ .
$$
Applying the Ricci identity for the tensors $T$ and $J$ in the above equality we obtain the assertion.
\end{proof}

In the rest of this section $x,y,z$ will be mutually orthogonal eigenvectors of $T$, which span a 
3-dimensional antiholomorphic plane in $T_pM$.
For any eigenvector $x$ of $T$ denote by $\lambda_x$  the corresponding eigenvalue.

\vspace{0.2cm}
{\bf Lemma 3.} {\it If for a triple $\{ x,y,z \}$
$$
	g((\nabla_{x}J)y,z)\ne 0 \ ,
$$
then the tensor $T$ is proportional to the metric tensor in the point $p$.}

\begin{proof} 
The substitution of $(u,v)$ by $(z,Jz)$ in (\ref{eq:2.6}) with the use of Lemma 1 and $M\in AK$ gives
\begin{equation} \label{eq:3.6}
	\begin{array}{l}
			\big(\nabla_x(T+Q)\big)(y,Jz)-\big(\nabla_y(T+Q)\big)(x,Jz)  \\
			+ \big( \ds \lambda_x-\frac12\lambda_y-\frac52\lambda_z+\mu\big) g((\nabla_xJ)y,z)
			+\big( \ds \frac12\lambda_x-\lambda_y+\frac52\lambda_z-\mu\big) g((\nabla_yJ)x,z)=0 \ .
	\end{array}
\end{equation}
Analogously from
$$
	(\nabla_{Jx}R)(y,x,x,z)+(\nabla_yR)(x,Jx,x,z)+(\nabla_xR)(Jx,y,x,z)=0	
$$
we find
\begin{equation} \label{eq:3.7} 
	\begin{array}{l}
			\big(\nabla_{Jx}(T+Q)\big)(y,z)-\big(\nabla_y(T+Q)\big)(Jx,z)  \\
			+ \big( \ds \lambda_x+\frac12\lambda_y+\frac12\lambda_z-\mu\big) g((\nabla_xJ)y,z)
			+\big( \ds -\frac92\lambda_x-\frac32\lambda_z+3\mu\big) g((\nabla_yJ)x,z)=0 \ .
	\end{array}
\end{equation}
From (\ref{eq:3.6}) and (\ref{eq:3.7}), using (\ref{eq:2.5}) and its analogue for $T$ we obtain
\begin{equation} \label{eq:3.8}
	(\lambda_y+3\lambda_z-2\mu)g((\nabla_xJ)y,z)+(-5\lambda_x+\lambda_y-4\lambda_z+4\mu)g((\nabla_yJ)x,z)=0 \ .
\end{equation}
Changing the places of $x$ and $y$ we have also
$$
	(\lambda_x-5\lambda_y-4\lambda_z+4\mu)g((\nabla_xJ)y,z)+(\lambda_x+3\lambda_z-2\mu)g((\nabla_yJ)x,z)=0 \ .
$$
Since $	g((\nabla_{x}J)y,z)\ne 0$, the last two equalities imply 
$$
	D(x,y,z)= 5 \lambda_x^2+ 5 \lambda_y^2- 7 \lambda_z^2 - 25 \lambda_x \lambda_y- 
 13 \lambda_x \lambda_z - 13 \lambda_y \lambda_z+ (14 \lambda_x  + 14 \lambda_y + 20 \lambda_z)\mu  -12 \mu^2=0 \ .
$$
Since $M$ is almost K\"ahler, at least one of $	g((\nabla_{y}J)z,x)$ and $	g((\nabla_{z}J)x,y)$ must also be
different from zero. Let e.g. $	g((\nabla_{y}J)z,x)\ne 0$. Then it follows $D(y,z,x)=0$. 

Case 1. $	g((\nabla_{z}J)x,y)$ also does not vanish. The $D(z,x,y)=0$. The system 
$$
	D(x,y,z)=D(y,z,x)=D(z,x,y)=0
$$
has a solution $x=y=z=\mu/2$.

Case 2. $	g((\nabla_{z}J)x,y)= 0$. Then (\ref{eq:3.8}) and $M\in AK$ imply
$$
	5\lambda_x-2\lambda_y+\lambda_z-2\mu=0 \ .
$$
The system 
$$
	D(x,y,z)=D(y,z,x)=5\lambda_x-2\lambda_y+\lambda_z-2\mu=0
$$
also has a solution $x=y=z=\mu/2$. If $n=3$, then $T=\mu/2g$ and the Lemma is proved.

Let $n\ge4$ and $u$ be an eigenvector of $T$, orthogonal to $span\{ x,y,z,Jx,Jy,Jz\}$. Replacing in (\ref{eq:2.6})
$(x,y,z,u,v)$ by $(y,z,u,Ju,x)$  we find 
$$ 
	(\lambda_x+\lambda_z+2\lambda_u-2\mu)g((\nabla_yJ)z,x)-(\lambda_x+\lambda_y+2\lambda_u-2\mu)g((\nabla_zJ)y,x)=0  
$$
and using $x=y=z=\mu/2$:
$$
	(2\lambda_u-\mu)g((\nabla_yJ)z,x)-(2\lambda_u-\mu)g((\nabla_zJ)y,x)=0  
$$
Now because of (\ref{eq:2.7}) it follows
$$
	(2\lambda_u-\mu)g((\nabla_xJ)y,z)=0  
$$
so $\lambda_u=\mu/2$, thus proving the Lemma.
\end{proof}

\vspace{0.2cm}
{\bf Lemma 4.} {\it If
$$
	g((\nabla_{y}J)y,x)\ne 0 \ ,
$$
then $\lambda_x=\lambda_z$.}

\begin{proof} 
We put in (\ref{eq:2.6}) $u=z$, $v=y$. Then we find
$$
	(\nabla_x(T+Q))(y,y)+(\nabla_x(T+Q))(z,z)=(\nabla_y(T+Q))(x,y)+(\nabla_z(T+Q))(x,z) \ .
$$
We replace here $z$ by $Jz$ and add the result with the above, using (\ref{eq:2.4}) and its
analogue for the tensor $T$. The result is
\begin{equation} \label{eq:3.9}
	(\nabla_x(T+Q))(y,y)=(\nabla_y(T+Q))(x,y) \ . 
\end{equation}
Now from the second Bianchi identity
$$
	(\nabla_{Jx}R)(x,y,y,x)+(\nabla_xR)(y,Jx,y,x)+(\nabla_yR)(Jx,x,y,x)=0	
$$
using (\ref{eq:3.9}), (\ref{eq:2.5}) and Lemma 1 we obtain
$$
	3\lambda_x+\lambda_y-2\mu=0  \ .
$$
Analogously we put in (\ref{eq:2.6}) $u=Jz$, $v=y$ and we derive
$$
	\lambda_x+\lambda_y+2\lambda_z-2\mu=0 \ .
$$
The last two equalities imply the assertion. 
\end{proof}

\vspace{0.2cm}
{\bf Lemma 5.} {\it Let $M$  be non K\"ahler in a point $p$, i.e. $(\nabla J)_p\ne o$. 
Then there exists a number $\lambda$, such that $T=\lambda g$ in $p$.}

\begin{proof}
Let  $\{ e_i,Je_i;\ i=1,...,n\}$ be an orthonormal basis of $T_pM$ of eigenvectors of $T$.
According to Lemmas  3 and 4 it suffices to consider the case
$$
	(\nabla_{e_1}J)e_1\ne 0 \qquad \qquad (\nabla_{e_i}J){e_j}= 0
$$
for any $i=2,...,n$, $j=1,...,n$. Moreover by Lemma 4 $\lambda_2=\lambda_3=...=\lambda_n$. 
Denote this number $\lambda$ and suppose $\lambda\ne\lambda_1$. Note that we may assume
$e_2\parallel (\nabla_{e_1}J)e_1$. Then we have also
$$
	g((\nabla_{e_1}J)e_1,e_i)=g((\nabla_{e_1}J)e_1,Je_i)= 0   \quad  {\rm for}\  i=3,...,n\ 
$$
and hence using again Lemma 3
\begin{equation} \label{eq:3.11}
	(\nabla_{e_1}J)e_i= 0 \quad  {\rm for}\  i=3,...,n\ .
\end{equation}
Putting in Lemma 2 \ $x=u=e_1$, $y=z=e_i$ $(i>1)$ we obtain
$$
	2\lambda_1+2\lambda+Q(e_1,e_1)+Q(e_i,e_i)-2\mu-\nu=0 \ . 
$$
Hence $Q(e_i,e_i)=Q(e_j,e_j)$ for any $i,j=2,...,n$. Making the same substitution
in (\ref{eq:2.3}) we find
$$
	Q(e_1,e_1)+Q(e_i,e_i)-\nu=-\frac12g((\nabla_{e_1}J)e_i,(\nabla_{e_1}J)e_i)
$$
for any $i=2,...,n$. Now $Q(e_2,e_2)=Q(e_3,e_3)$ and (\ref{eq:3.11}) imply
$$
	g((\nabla_{e_1}J)e_2,(\nabla_{e_1}J)e_2)=0
$$
and hence
$$
	(\nabla_{e_1}J)e_1=0 \ ,
$$
which is a contradiction. This proves the lemma.
\end{proof}

\vspace{0.2cm}
{\bf Theorem 2.} {\it Let $M$ be a $2n$-dimensional conformal flat $AK_3$-manifold, $n\ge 2$. 
Then $M$ is a flat K\"ahler manifold or $n=2$ and $M$ is locally a product
of two 2-dimensional K\"ahler manifolds $M_1$ and $M_2$ of constant sectional 
curvature $c$ and $-c$, $c>0$, respectively. }

\begin{proof} 
According to \cite{OK2} and \cite{OK3} a  conformal flat $AK_3$-manifold 
of dimension $\ge 4$ must be a 4- or a 6-dimensional manifold of
constant sectional curvature, or a flat K\"ahler manifold, or a product
of two almost K\"ahler manifolds $M_1$ and $M_2$ of constant sectional 
curvature $c$ and $-c$, $c>0$, respectively. On  the other hand by \cite{O-S}
an almost K\"ahler manifold of constant sectional curvature and dimension
$\ge4$ is a flat K\"ahler manifold. Hence the assertion follows.
\end{proof}

{\bf Remark.} There exist conformal flat almost K\"ahler manifolds which
are not K\"ahler, see e.g. \cite{Catalano--Olszak}. So the $AH_3$-assumption
can not be removed.

\vspace{0.2cm}
{\bf Theorem 3.} {\it Let $M$ be a $2n$-dimensional $AK_3$-manifold, $n\ge 3$,
with vanishing Tricerri-Vanhecke Bochner curvature tensor and constant scalar curvature. 
Then $M$ is a K\"ahler manifold. }

\begin{proof} 
Assume that $M$ is non K\"ahler in a point $q$. Then $M$ is 
non K\"ahler in a neighborhood $U$ of $q$. By  Lemma 5  in $U$ 
it holds $T=\lambda g$ with a function $\lambda$. Now from
$B=0$ it follows that in $U$ the curvature tensor of $M$ has the form
$$
	R=\varphi(Q)+\theta(\pi_1+\pi_2)-\nu\pi_1 
$$
with a function $\theta$. According to Theorem 1 $U$ must be conformal flat. 
By Theorem 2 this contradicts the assumption that $U$ is non K\"ahler. 
\end{proof}

\end{document}